\documentclass[a4paper,twoside]{article}

\def\shorttitle{}
\def\shortauthor{Dexie Lin}

\usepackage{amsfonts}
\usepackage{amsmath}
\usepackage{amssymb}
\usepackage[all,cmtip]{xy}
\usepackage{amscd}
\usepackage{amsbsy}
\usepackage{amsthm}
\usepackage{bm}
\usepackage{empheq}
\usepackage{graphicx}
\usepackage{psfrag}
\usepackage{color}
\usepackage{cite}
\usepackage{slashed}
\usepackage{multicol}
\usepackage[colorlinks,linkcolor=blue,citecolor=red]{hyperref}
\usepackage{xcolor,graphicx,blkarray}
\usepackage{hyperref}

\allowdisplaybreaks[4]

\newfont{\myfnt}{cmssi10 scaled 1440}
\numberwithin{equation}{section}
\catcode`@=11
\def\ps@nk{\def\@oddhead{\vbox{\hbox to \hsize{\pic \footnotesize \it \shorttitle
\hfill \rm \thepage} \vspace{1mm} \vspace*{-2mm}}}
\def\@evenhead{\vbox{\hbox to \hsize{\pic \footnotesize \rm \thepage \hfill \it \shortauthor}
\vspace{1mm} \vspace*{-2mm}}}
\def\@oddfoot{} \def\@evenfoot{}}
\def\ps@first{\def\@oddhead{\vbox{\hbox to \hsize{\pic \footnotesize
} \break}}
\def\@oddfoot{} \def\@evenfoot{}}
\newtheoremstyle{thmstyle}
  {6pt}
  {6pt}
  {\it}
  {}
  {\bf}
  {}
  {.5em}
  {}
\newtheoremstyle{remstyle}
  {6pt}
  {6pt}
  {\rm}
  {}
  {\bf}
  {}
  {.5em}
  {}

\def\Sec{\@Startsection{section}{1}{\z@}
                                   {-3.5ex \@plus -1ex \@minus -.2ex}%
                                   {2.3ex \@plus.2ex}%
                                   {\normalfont\large\bfseries\boldmath}}
\def\@Startsection#1#2#3#4#5#6{%
  \if@noskipsec \leavevmode \fi
  \par
  \@tempskipa #4\relax
  \@afterindenttrue
  \ifdim \@tempskipa <\z@
    \@tempskipa -\@tempskipa \@afterindentfalse
  \fi
  \if@nobreak
    \everypar{}%
  \else
    \addpenalty\@secpenalty\addvspace\@tempskipa
  \fi
  \@ifstar
    {\@ssect{#3}{#4}{#5}{#6}}%
    {\@dblarg{\@Sect{#1}{#2}{#3}{#4}{#5}{#6}}}}
\def\@Sect#1#2#3#4#5#6[#7]#8{%
  \ifnum #2>\c@secnumdepth
    \let\@svsec\@empty
  \else
    \refstepcounter{#1}%
    \protected@edef\@svsec{\@seccntformat{#1}\relax}%
  \fi
  \@tempskipa #5\relax
  \ifdim \@tempskipa>\z@
    \begingroup
      #6{%
          \@hangfrom{\hskip #3\relax\@svsec \hskip -2.5mm}%
          \interlinepenalty \@M #8\@@par}
    \endgroup
    \csname #1mark\endcsname{#7}%
    \addcontentsline{toc}{#1}{%
      \ifnum #2>\c@secnumdepth \else
        \protect\numberline{\csname the#1\endcsname}%
      \fi
      #7}%
  \else
    \def\@svsechd{%
      #6{\hskip #3\relax
      \@svsec #8}%
      \csname #1mark\endcsname{#7}%
      \addcontentsline{toc}{#1}{%
        \ifnum #2>\c@secnumdepth \else
          \protect\numberline{\csname the#1\endcsname}%
        \fi
        #7}}%
  \fi
  \@xsect{#5}}
\renewenvironment{abstract}{%
        \small
        \quotation
         \noindent {\bfseries \abstractname } }%
      {\if@twocolumn\else\endquotation\fi}

\def\Subsec{\@StartSubsection{subsection}{2}{\z@}%
                                     {-3.25ex\@plus -1ex \@minus -.2ex}%
                                     {1.5ex \@plus .2ex}%
                                     {\normalfont\normalsize\bfseries\boldmath}}
\def\@StartSubsection#1#2#3#4#5#6{%
  \if@noskipsec \leavevmode \fi
  \par
  \@tempskipa #4\relax
  \@afterindenttrue
  \ifdim \@tempskipa <\z@
    \@tempskipa -\@tempskipa \@afterindentfalse
  \fi
  \if@nobreak
    \everypar{}%
  \else
    \addpenalty\@secpenalty\addvspace\@tempskipa
  \fi
  \@ifstar
    {\@ssect{#3}{#4}{#5}{#6}}%
    {\@dblarg{\@SubSect{#1}{#2}{#3}{#4}{#5}{#6}}}}
\def\@SubSect#1#2#3#4#5#6[#7]#8{%
  \ifnum #2>\c@secnumdepth
    \let\@svsec\@empty
  \else
    \refstepcounter{#1}%
    \protected@edef\@svsec{\@seccntformat{#1}\relax}%
  \fi
  \@tempskipa #5\relax
  \ifdim \@tempskipa>\z@
    \begingroup
      #6{%
          \@hangfrom{\hskip #3\relax\@svsec\hskip -1.5mm}%
          \interlinepenalty \@M #8\@@par}
    \endgroup
    \csname #1mark\endcsname{#7}%
    \addcontentsline{toc}{#1}{%
      \ifnum #2>\c@secnumdepth \else
        \protect\numberline{\csname the#1\endcsname}%
      \fi
      #7}%
  \else
    \def\@svsechd{%
      #6{\hskip #3\relax
      \@svsec #8}%
      \csname #1mark\endcsname{#7}%
      \addcontentsline{toc}{#1}{%
        \ifnum #2>\c@secnumdepth \else
          \protect\numberline{\csname the#1\endcsname}%
        \fi
        #7}}%
  \fi
  \@xsect{#5}}
\def\list#1#2{\ifnum \@listdepth >5\relax \@toodeep \else \global
\advance \@listdepth\@ne \fi \rightmargin \z@ \listparindent\z@
\itemindent\z@ \csname @list\romannumeral\the\@listdepth\endcsname
\def\@itemlabel{#1}\let\makelabel\@mklab \@nmbrlistfalse #2\relax
\@trivlist \parskip 0pt \parindent\listparindent \advance \linewidth
-\rightmargin \advance\linewidth -\leftmargin \advance\@totalleftmargin
\leftmargin \parshape \@ne \@totalleftmargin \linewidth \ignorespaces}
\renewcommand{\@makecaption}[2]{\begin{center}#1. #2\end{center}}
\catcode`@=12 
\theoremstyle{thmstyle}
\newtheorem{thm}{\indent Theorem}[section]
\newtheorem{lemma}[thm]{\indent Lemma}
\newtheorem{prop}[thm]{\indent Proposition}

\theoremstyle{remstyle}

\def\Dirac{\slashed{D}}

\newsavebox{\mygraphic}
\def\pic{\begin{picture}(0,0) \put(-210,-1250){\usebox{\mygraphic}} \end{picture}}
\newfont{\HUGEbf}{cmbx10 scaled 3500}
\definecolor{gray}{rgb}{0.9,0.9,0.9}
\def\thebibliography#1{\section*{\bf \large References}
\list{[\arabic{enumi}]} {\settowidth \labelwidth{[#1]} \leftmargin
\labelwidth \advance \leftmargin \labelsep \usecounter{enumi}}
\def\newblock{\hskip .11em plus .33em minus .07em} \footnotesize \sloppy \clubpenalty
4000 \widowpenalty 4000 \sfcode`\.=1000 \relax}

\def\spinc{$\mbox{spin}^c$}

\theoremstyle{definition}

\numberwithin{equation}{section}

\title{Seiberg-Witten equation on a manifold with rank  $2$-foliation}

\author{\large Dexie Lin}

\date{}

\begin{document}

\maketitle
\date{}
\thispagestyle{first}
\renewcommand{\thefootnote}{\fnsymbol{footnote}}

\footnotetext{\hspace*{-5mm} \begin{tabular}{@{}r@{}p{14cm}@{}} &
Manuscript last updated: \today.\\

\end{tabular}}

\renewcommand{\thefootnote}{\arabic{footnote}}

\begin{abstract}
Let $M$ be a closed oriented $4$-manifold admitting a rank-$2$ oriented foliation with a metric of leafwise positive scalar curvature. If $b^+>1$, then we will show that   the Seiberg-Witten invariant vanishes for all  \spinc structures.

\vskip 4.5mm

\noindent\begin{tabular}{@{}l@{ }p{10cm}} {\bf Keywords } & Foliation, Bundle-like metric, Seiberg-Witten invariant
\end{tabular}

\vskip 4.5mm


\end{abstract}

\baselineskip 14pt

\setlength{\parindent}{1.5em}

\setcounter{section}{0}

\section{Introduction}



A natural question in Riemannian geometry is: When does a closed manifold $M$ admit
a Riemannian metric with positive scalar curvature? From the  Lichnerowicz formula \cite{L}, we know that $\hat A(M)$ is  a topological obstruction for a manifold to admit a metric of positive scalar curvature. A manifold equipped with a foliation is called a foliated manifold. We can define the leafwise scalar curvature as follows: For any $x\in M$, the integrable subbundle $F_x$ determines a leaf $\mathcal F_x$, such that $T\mathcal F_x=F_x$. Hence the metric on $M$ determines a Riemannian metric on $\mathcal F_x$. Using the Levi-Civita connection, one has the scalar  curvature of this metric, which is denoted by  $scal^{\mathcal F_x}$.

 By the  results given by A. Connes \cite{Connes} and W. Zhang \cite{Zhang} respectively, it is known  that
 $\hat A(M)$ is also a topological obstruction  for the foliation admitting metric of positive leafwise scalar curvature.
\begin{thm}[A. Connes 1986]
  For a closed oriented manifold, let $\mathcal F$ be a spin foliation. If $\hat A(M)\neq0$, then $\mathcal F$ does not admit any metric with positive scalar curvature.
\end{thm}
\begin{thm}[W. Zhang 2017]
  For a closed oriented spin manifold, let $\mathcal F$ be a   foliation. If $\hat A(M)\neq0$, then $\mathcal F$ does not admit any metric with positive scalar curvature.
\end{thm}

\noindent
Later Zhang \cite{Zhang2} posed the following question.

\vspace{3mm}

 {\bf Question:} Given a closed oriented 4-manifold with a foliation with positive leafwwise scalar
curvature, does the Seiberg-Witten invariant vanish?

\vspace{3mm}
\noindent
On the other hand,
for oriented smooth $4$-manifolds,  Seiberg-Witten invariant plays an important role to study the obstruction for manifold admitting some geometric and topological structures. One well-known result is that   a closed oriented $4$-manifold with non-trivial Seiberg-Witten invariant  can not admit any metric of positive scalar curvature.

\noindent
 In this paper we will give a partial answer to the above question.
Let $M$ denote  a  closed oriented smooth $4$-manifold  with $b^+>1$ and admit a rank-$2$ foliation $\mathcal F$, where $b^+$ denotes the dimension of the self-dual part of the second cohomology groups, see \cite[Section 1]{KM}. Suppose that $\mathcal F$ admits a positive scalar curvature metric. Then, we show that each leaf is compact  and   the following  result holds.

\begin{thm}\label{thm-vanishing-1}
  Let $M$ denote  a  closed oriented smooth $4$-manifold  with a rank-$2$  foliation admitting a leafwise positive scalar curvature metric. Suppose that the positive part of the second Betti number is greater than $1$, i.e.  $b^+>1$. Then the Seiberg-Witten invariant vanishes for all \spinc structures.
\end{thm}



The structure of this paper is as follows: in Section 2, we will show that each leaf is compact and the manifold $M$ admits a bundle-like metric under the condition of Proposition \ref{prop-trivial-1}; in Section 3, we review the classical theory of the Seiberg-Witten equation and  give a proof to Theorem \ref{thm-vanishing-1} via the adiabatic limit method.

\vspace{5mm}

{\bf Acknowledgement} The author would like to express the special thanks to Huitao Feng for the introduction of such a problem and Mikio Furuta for the invaluable discussion. The author wants to thank Clifford Taubes for sharing  the idea to proving the compactness of the leaf.   The research is supported by the FMSP.

\section{Compact leaf and bundle-like metric}

\subsection{Compact leaf}

In this subsection, we   show that each leaf in the manifold satisfying the condition of Theorem \ref{thm-vanishing-1} must be compact.  In general, we can establish the following proposition.

\begin{prop}\label{prop-compactness}
  Let $M$ denote  a  closed   smooth $n$-manifold with  a rank-$2$   foliation.  Suppose that  the foliation admits a metric with leafwise positive scalar
curvature. Then, each leaf is compact. Moreover, each leaf is diffeomorphic  to either $S^2$ or $ \mathbb RP^2$.
\end{prop}

\begin{pf}
From the compactness of $M$, we know that each leaf is a complete Riemannian submanifold with the induced metric.
It suffices to show that for a given geodesic in a leaf, there always exists a conjugate point in this geodesic.
  Since the dimension of each leaf is $2$,  the leafwise sectional curvature equals half of the leafwise scalar curvature,  which implies that the sectional curvature of each leaf stays away from zero, i.e. there exists some number $c_0>0$ such that over each leaf the leafwise sectional curvature is strictly greater than $c_0$.  Let $\mathcal F_x$ denote the leaf along a point $x\in M$; we need to show that  it is compact. Choose a geodesic $\gamma$ starting at $x$ in  $\mathcal F_x$ and consider the Jacobi field along this geodesic
  \[\frac{d^2J}{dt^2}+\kappa J=0,~J(0)=0,\]
  where $\kappa$ denotes the sectional curvature.
  By \cite[Chapter 5]{Carmo}, we can write the Jacobi field as follows:
  $$J(t)=f(t)e(t),$$
   where
  $e(t)$ is the parallel normal vector field along the geodesic, i.e. $\gamma'(t)\perp e(t)$ and $\frac{de(t)}{dt}=0$ and $f$ is a function along the geodesic such that $f(0)=0$.   We rewrite the Jacobi equation as
  \[\frac{d^2f(t)}{dt^2}+\kappa(t)f(t)=0.\]
 We claim that $f(t)$ must vanish at some point of the geodesic. Suppose it is not true, without loss of generality we assume that $f(t)>0$. Consider the derivative of $f$, we can categorize it into two cases:
  \begin{enumerate}
    \item There exists $t_0>0$ such that $f'(t_0)<0$.  By the equation $f''(t)=(f'(t))'=
  -\kappa(t) f(t)<0$, we have $f'(t)< f'(t_0)<0$ for any $t\geq t_0$. Thus, there exists some point $t_1$ such that $f(t_1)<0$, which contradicts to the hypothesis.
    \item For all $t$, we have  $f'(t)\geq0$. By choosing a small enough $t'_0>0$, we have $f''(t)=-\kappa(t)f(t)\leq -c f(t'_0)$ for $t\geq t'_0$. Apply the similar argument, we can   deduce that there exists a point $t_2$ such that $f'(t_2)<0$, which contradicts to the hypothesis.
  \end{enumerate}
  By the above argument, we can always find a point such that $f$ vanishes.
  Hence,  each leaf is compact. By the classification of closed $2$-manifold with positive sectional curvature, one has that each leaf is diffeomorphic to either $S^2$ for the orientable leaf or $\mathbb RP^2$ for the non-orientable leaf.
\end{pf}

Under the condition of Theorem \ref{thm-vanishing-1}, we will show that for each leaf there is a saturated neighborhood tubular neighborhood as defined in \cite{Thur}. In particular,
for the oriented foliation case, we can show that the manifold becomes a fibration if the foliation admits a metric with  positive scalar
curvature.

\begin{prop}\label{prop-trivial-1}
  Let $M$ denote  a  closed   smooth $n$-manifold with  a rank-$2$ oriented   foliation $\mathcal F$.  Suppose that  the foliation admits a metric with leafwise positive scalar
curvature. Then, $M$ is a   fibration over some closed manifold $B$ whose fiber is diffeomorphic to $S^2$ i.e. there exists a submersion
$$\pi: M{\to} B$$   and for each point $p\in B$ we have a diffeomorphism $\pi^{-1}(p)\cong S^2$.
\end{prop}

\begin{pf} It suffices to show that $M/\mathcal F$ is a smooth manifold.
  Proposition \ref{prop-compactness} tells that each leaf is $S^2$. Fixing a leaf $\mathcal L$  and local transversal
$T$, we have a leaf holonomy \cite[Chapter 1.7]{Molino}
 $$ Hol:\pi_1(\mathcal L,x_0)\to Diff_{x_0}(T),$$
   from the fundamental group of this leaf for the fixed point $x_0\in \mathcal L$ to the germs of the local diffeomorphism of the transverse manifold at $x_0$.  Since $\pi_1(S^2,x_0)$ is trivial, this holonomy is trivial. Therefore, on each leaf $\mathcal L$, there is a  neighborhood $N(\mathcal L)$ which is diffeomorphic to the standard product $S^2\times D^2$. It implies that the quotient  $M/\mathcal F$ has a smooth manifold structure.
\end{pf}

In general, for the non-orientable foliation case, we have the following proposition.

\begin{prop}\label{prop-trivial-2}
   Let $M$ denote  a  closed oriented smooth $4$-manifold  with a rank-$2$  foliation $\mathcal F$.  Suppose that  the foliation admits a metric with leafwise positive scalar
curvature.  Then $M/\mathcal F$ is an orbifold.
\end{prop}

\begin{pf}
  We give a sketch of the proof. If each leaf is orientable, then it is known that each leaf is diffeomorphic to  $S^2$, and the above arguments of  Proposition \ref{prop-trivial-1} assert that $M/\mathcal F$ is a manifold. In general, we consider the double-covering $\tilde M$ of $M$ with the canonical lift foliation $\tilde{\mathcal F}$. We have that $\tilde M/\tilde{\mathcal F}$ is a manifold, and the double covering
  \[\tilde M/\tilde{\mathcal F}\to M/\mathcal F\] equips an orbifold structure for $M/\mathcal F$.
\end{pf}

Combining Proposition \ref{prop-trivial-1} and Proposition \ref{prop-trivial-2}, we have the following proposition.

\begin{prop}\label{prop-bundle-like}
  Let $M$ denote  a  closed oriented smooth $4$-manifold  with a rank-$2$  foliation $\mathcal F$. Suppose that  the foliation admits a metric with leafwise positive scalar
curvature.   Then $M$ admits a bundle-like metric such that the restriction to the foliation coincides with $g_{\mathcal F}$.
\end{prop}

\begin{pf}
  By Proposition \ref{prop-trivial-2}, we have that $M/\mathcal F$ is an orbifold.
   It is well-known that any orbifold admits a Riemannian metric. We can pull back a metric of $M/\mathcal F$ to an $\mathcal F$-transverse metric on $M$ and  complete to a bundle-like metric.
\end{pf}


{\bf Remark}: \begin{enumerate}
   \item Following the same idea, one can show that:  If a closed manifold $M$ with foliation $\mathcal F$ has each compact leaf  and  has finite leafwise holonomy, then $M/\mathcal F$ is an orbifold.
   \item Conversely, if each leaf is compact and $M$ admits a bundle-like metric, then $M/\mathcal F$ is an orbifold(c.f. \cite[Proposition 3.7]{Molino}).
   \item  If the holonomy is not finite, there is an example(c.f. \cite{Sull}) such that even though the leaves are compact,   the quotient space $M/\mathcal F$ fails to be an orbifold.
 \end{enumerate}

\subsection{Bundle-like metrics }
For a foliated manifold, a notion of bundle-like metric was firstly posted by  Reinhart \cite{Reinhart}.
 Let $\mathcal F$ be an integrable subbundle of the tangent vector bundle $TM$ of a smooth Riemannian
manifold $(M,g)$. Then, we have the associated foliation $\mathcal F$. There is a splitting for this metric
\[g=g_{\mathcal F}\oplus g_{\mathcal F^\perp},\]
and an isomorphism
\[\mathcal F^\perp\cong Q,\]
where $Q$ denotes the quotient $TM/ \mathcal F$. $Q$ inherits a metric $g_Q=g_{ \mathcal F^\perp}$. We say $g_Q$ is \emph{bundle-like}, if
\[L_vg_Q\equiv0,~\mbox{for all }v\in \Gamma(\mathcal F),\]
here $L_v$ denotes the Lie-derivative associated with $v$.
Given a bundle-like metric,
we define \[g_\beta=(\beta^2g_{\mathcal F})\oplus g_Q.\]
Denote by $\nabla^\beta$ the associated  Levi-Civita connection and $\langle \rangle$ the metric of $g_0$.

\begin{lemma}
 By the straightforward calculation, we have that \begin{itemize}
     \item[(1)] \[\langle\nabla^{\beta}_{e_i}e_j,e_k\rangle=O(1),~
     \langle\nabla^{\beta}_{e_i}e_j,f_k\rangle=O(\beta^2).\]
     \item[(2)] \[\langle\nabla^{\beta}_{e_i}f_j,e_k\rangle=O(1),~
     \langle\nabla^{\beta}_{e_i}f_j,f_k\rangle=O(1).\]
     \item[(3)]
     \[\langle\nabla^{\beta}_{f_i}e_k,e_j\rangle=O(1),~
     \langle\nabla^{\beta}_{f_i}e_j,f_k\rangle=O(\beta^2).\]
      \item[(4)]
     \[\langle\nabla^{\beta}_{f_i}f_k,e_j\rangle=O(1),~
     \langle\nabla^{\beta}_{f_i}f_j,f_k\rangle=O(1).\]
 \end{itemize}
 where  $e_i\in\Gamma(F)$, $f_i\in\Gamma(Q)$.
\end{lemma}

\begin{thm}\label{thm-scalar-curvature}
The scalar curvature $Scal^\beta$ associated with the metric $g_\beta$ can be expressed as follows:
\[Scal^\beta=\frac{Scal^F}{\beta^2}+O(1).\]
\end{thm}
\begin{pf} Let $p: TM\to F$ and $p^\perp: TM\to F^\perp$ be the orthogonal projection maps.
For $e_i,e_j\in\Gamma(F)$, we get
  \begin{eqnarray*}
    \langle R^\beta(e_i,e_j)e_i,e_j\rangle&=&
    \langle\nabla^\beta_{e_i}(p+p^\perp)\nabla^\beta_{e_j}e_i,e_j\rangle
    -\langle\nabla^\beta_{e_j}(p+p^\perp)\nabla^\beta_{e_i}e_i,e_j\rangle-
    \langle \nabla^\beta_{[e_i,e_j]}e_i,e_j\rangle\\
    &=&\langle R^F(e_i,e_j)e_i,e_j\rangle-\beta^2\langle p^\perp\nabla_{e_j}e_i,\nabla_{e_i}e_j\rangle+\beta^2\langle
    p^\perp\nabla_{e_i}e_i,\nabla_{e_j}e_j\rangle\\
    &=&\langle R^F(e_i,e_j)e_i,e_j\rangle+O(\beta^2).
  \end{eqnarray*}
  For $e_i\in\Gamma(F),~f_j\in\Gamma(Q)$, we have
  \begin{eqnarray*}
    \langle R^\beta(e_i,f_j)e_i,f_j\rangle&=&\beta^2\langle\nabla _{e_i}p\nabla_{f_j}e_i,f_j\rangle+\beta^2\langle\nabla _{e_i}p^\perp\nabla_{f_j}e_i,f_j\rangle-\beta^2\langle\nabla _{f_j}p\nabla_{e_i}e_i,f_j\rangle\\
    &&-\beta^2\langle\nabla _{f_j}p^\perp\nabla_{e_i}e_i,f_j\rangle-
    \beta^2\langle\nabla_{[e_i,f_j]}e_i,f_j\rangle\\
    &=&O(\beta^2).
  \end{eqnarray*}
  Similarly, for $f_i,~f_j\in \Gamma(Q)$, we have
  \begin{eqnarray*}
    \langle R^\beta(f_i,f_j)f_i,f_j\rangle&=&
    \beta^2\langle\nabla_{f_i}p\nabla^\beta_{f_j}f_i,f_j\rangle+
    \langle\nabla_{f_i}p^\perp\nabla_{f_j}f_i,f_j\rangle-\beta^2
    \langle\nabla_{f_j}p\nabla^\beta_{f_i}f_i,f_j\rangle\\
    &&-\langle\nabla_{f_j}p^\perp\nabla_{f_i}f_i,f_j\rangle-
    \langle\nabla_{[f_i,f_j]}f_i,f_j\rangle\\
    &=&O(1).
  \end{eqnarray*}
  Combining the above three formulas, one gets the desired result.
\end{pf}

\section{Seiberg-Witten invariant and  vanishing theorem }
In this section, we will   give a proof of Theorem \ref{thm-vanishing-1}. The idea to show the theorem is as follows: By Proposition \ref{prop-bundle-like}, if the hypothesis of Theorem \ref{thm-vanishing-1} is satisfied, then in Theorem \ref{thm-scalar-curvature} we may choose $\beta$ small enough so that the scalar curvature of the manifold is positive. Thus, the   conclusion of
Theorem \ref{thm-vanishing-1} follows.
Before proceeding, we review the Seiberg-Witten invariant and the result that if a closed  oriented $4$-manifold has positive scalar curvature and $b^+>1$, then the Seiberg-Witten invariant vanishes.

Let $(M,g)$ be a closed oriented Riemannian $4$-manifold with a \spinc structure $\mathfrak s$. Let $S^\pm$ denote the spinor bundles associate to $\mathfrak  s$, there is a well-defined Dirac operator
\[\Dirac_A:\Gamma(S^+)\to\Gamma(S^-),\]
where $A$ is a connection on the determinant line bundle of this \spinc structure $\mathfrak s$.
We give a brief introduction to the classical Seiberg-Witten theory, see \cite[Chaper 3,4]{Morgan} for more details. Let $\mathcal A$ denote all the connections on the determinant line bundle,
for $(A,\Phi)\in\mathcal A\times \Gamma(S^+)$, we define Seiberg-Witten equation as follows:
\[\begin{cases}
  F^+_A=q(\Phi)\\
  \Dirac_A\Phi=0\\
\end{cases},\]
where $q(\Phi)=\Phi\otimes\Phi^*-\frac{|\Phi|^2}21$ and we used the identification
\[cl_+:\Lambda^{2,+}\otimes\mathbb C\to End^0(S^+),\]
between the self-adjoint two forms and traceless  endomorphism  of $S^+$.
The moduli space $\mathcal M(\mathfrak s)$ is the space consisting  of all solutions $(A,\Phi)$ mod  the   gauge group $\mathcal G=C^\infty(M,S^1)$, where the gauge action is defined by the following: for each $g\in\mathcal G$,
\[g(A,\Phi)=(A-2g^{-1}dg,g\cdot \Phi),\]
 where $\cdot$ denotes the Clifford multiplication. We can also perturb the equations, by adding a self-dual two-form $\eta$, namely to solve the equations
\[\begin{cases}
  F^+_A=q(\Phi)+\eta\\
  \Dirac_A\Phi=0\\
\end{cases},\]
we write $\mathcal M_\eta(\mathfrak s)$ for the perturbed moduli space. The formal dimension of the moduli space $\mathcal M_\eta(\mathfrak s)$ is
\[d(\mathfrak s)=\frac14(c_1(\det(\mathfrak s))\cdot c_1(\det(\mathfrak s))-2e(M)
-3\sigma(M)),\]
where $e(M)$ denotes the Euler number and $\sigma(M)$ denotes the signature.

We have the following well-known results about the moduli space: \begin{enumerate}
  \item The moduli space $\mathcal M_\eta(\mathfrak s)$ is compact.
  \item The orientation of  $H^0(M)\otimes H^1(M)\otimes H^+(M)$ induces the orientation of the moduli space $\mathcal M_\eta(\mathfrak s)$
  \item For an open dense set of the perturbation $\eta$, the moduli space is a smooth manifold consisting of irreducible solution(i.e.  $\Phi\neq0$). 
\end{enumerate}
By fixing an orientation of the moduli space,  we define the Seiberg-Witten invariant to be zero, if the formal dimension is odd or less than zero, otherwise the Seiberg-Witten invariant $\mathfrak s$ is to be the following:
\[SW_\eta(\mathfrak s)=\int_{\mathcal M_\eta(\mathfrak s)}\mu^{d/2},\]
where $\mu$ denotes the first Chern class of the canonically associated principal $S^1$-bundle, i.e. the solutions mod the reference gauge group $\mathcal G_0=\{u\in \mathcal G\big| u(x_0)\equiv1\}$, where $x_0\in M$ is a fixed point.
If $b^+(M)>1$, then it is known that the moduli space is generically  independent of the  choice of the perturbation and the metric. In this case, we often omit the subscript of the perturbation.

We review the classical result of the local estimate for Seiberg-Witten equation \cite[Chapter 4]{Morgan}.
Let $(A,\Psi)$ be solution of the Seiberg-Witten equation, by Weitzenb\"ock formula one deduces that
\[
0=\frac12\Delta|\Psi|^2+|\nabla_A\Psi|^2+\frac14Scal|\Psi|^2+\frac14|\Psi|^4,
\]
hence\[\Delta|\Psi|^2+\frac12Scal|\Psi|^2\leq0.\]
 At the maximal point of $|\Psi|^2$, one gets that $\Psi=0$ in order that the
scalar curvature is positive, so that the Seiberg-Witten invariant vanishes.
Together with Proposition \ref{prop-bundle-like}, Theorem \ref{thm-scalar-curvature}, and letting $\beta\to0$, we  proved Theorem \ref{thm-vanishing-1}. Therefore, the Seiberg Witten invariant is an obstruction for the closed oriented $4$ manifold admitting a rank $2$-foliation with positive scalar curvature.

Graduate School of Mathematical Sciences, The University of Tokyo, 3-8-1 Komaba, Meguro-ku, Tokyo 153-8914, Japan. E-mail: \mbox{dexielin@ms.u-tokyo.ac.jp}
\end{document}